# EXACT SEQUENCES OF FIBRATIONS OF CROSSED COMPLEXES, HOMOTOPY CLASSIFICATION OF MAPS, AND NONABELIAN EXTENSIONS OF GROUPS

RONALD BROWN

*Abstract*

The classifying space of a crossed complex generalises the construction of Eilenberg-Mac Lane spaces. We show how the theory of fibrations of crossed complexes allows the analysis of homotopy classes of maps from a free crossed complex to such a classifying space. This gives results on the homotopy classification of maps from a $CW$-complex to the classifying space of a crossed module and also, more generally, of a crossed complex whose homotopy groups vanish in dimensions between 1 and $n$. The results are analogous to those for the obstruction to an abstract kernel in group extension theory.

## Introduction

A major result in the theory of crossed complexes is the homotopy classification theorem of [**BH91**]:

$$[\Pi X_*, C] \cong [X, BC] \tag{1}$$

where:

- on the left, $\Pi X_*$ is the fundamental crossed complex of the skeletal filtration of the CW-complex $X$, and square brackets denote homotopy classes in the category Crs of crossed complexes, and

- on the right, $BC$ is the classifying space of the crossed complex $C$, and square brackets denote homotopy classes of maps of topological spaces.

Results of Eilenberg–Mac Lane, including the case when $C$ is just a group, and also the case of a local coefficient system, are obtained by giving special cases for the crossed complex $C$. A new application is to the case $C$ is a crossed module, and this is exploited in, for example, [**FK08, Far08, PT07**].

The purpose of this paper is to show how this homotopy classification can be analysed using the notion of exact sequence of a fibration of crossed complexes. In this way we show how crossed complexes give a convenient common setting for the homotopy classification of maps to an $n$-aspherical space $Y$ (i.e. one in which $Y$ is path connected and $\pi_i Y = 0$ for $1 < i < n$) as originally found by Olum in [**Olu53**], and for the classical theory of nonabelian extensions of groups and abstract kernels, [**ML63**]. A topological account of this extension theory is given in [**Ber85**].

An example of the convenience of the setting of crossed complexes in these areas is that the traditional Schreier notion of factor set for an extension $1 \to K \to E \to G \to 1$ of groups is conveniently formulated as a morphism of crossed complexes $F^{st}(G) \to \text{AUT}(K)$ from the standard free crossed resolution $F^{st}(G)$ of $G$ to the crossed module $\text{AUT}(K) := (K \to \text{Aut}(K))$, and equivalence of factor sets is just homotopy of morphisms, [**BH82**]. By replacing $F^{st}(G)$ by an equivalent free crossed resolution, we can get calculations, as exploited in [**BP96**]. That exposition also shows the use of extensions of the type of a crossed module, due originally to Dedecker, [**Ded58**]. This relates to work of [**BM94**].

Our techniques involve the closed monoidal structure on the category of crossed complexes from [**BH87**], the exact sequences of a fibration of crossed complexes established by Howie in [**How79**], and the model category properties of crossed complexes from [**BG89, BH91**]. For more on model categories, see [**Hov99**].

This paper should also be seen as a modern development of pioneering work of J.H.C. Whitehead in [**Whi49**], which in some ways goes further than the work of Olum published later. Whitehead uses the term 'homotopy system' for what we call 'reduced, free crossed complex'. These ideas are exploited by Ellis in [**Ell88**]. Although Whitehead does not have the notion of classifying space, he gives his view that, in our terminology, crossed complexes have better realisation properties than chain complexes with a group of operators.

We are thus confirming that crossed complexes, though only a linear model of homotopy types, are convenient in that border between homology and homotopy which includes the Relative Hurewicz Theorem, [**BH81**], the Homotopy Addition Lemma [**BS07**], computation with crossed modules [**BW95**], and this level of the homotopy classification of maps.

The referee is thanked for some helpful comments.

## 1.  Crossed complexes

A crossed complex $C$ is in part a sequence of the form

$$\cdots \longrightarrow C_n \xrightarrow{\delta_n} C_{n-1} \xrightarrow{\delta_{n-1}} \cdots \cdots \xrightarrow{\delta_3} C_2 \xrightarrow{\delta_2} C_1$$

where all the $C_n, n \geqslant 1$ are groupoids over $C_0$. The structure and axioms for a crossed complex are those universally satisfied by the main topological example, the *fundamental crossed complex* $\Pi X_*$ of a filtered space $X_*$, where $(\Pi X_*)_1$ is the fundamental groupoid $\pi_1(X_1, X_0)$ and $(\Pi X_*)_n$ is the family of relative homotopy groups $\pi_n(X_n, X_{n-1}, x_0)$ for all $x_0 \in X_0$.

The fundamental groupoid of the crossed complex $C$ is $\pi_1 C = \text{Coker}\, \delta_2$ and the homology groups of $C$ are for $n \geqslant 2$ the families of abelian groups

$$H_n(C, x) = (\text{Ker}\, \delta_n(C_n(x) \to C_{n-1}(x))) / (\delta_{n+1} C_{n+1}(x)) \qquad x \in C_0.$$

The crossed complex $C$ is *aspherical* if it is connected (i.e. the groupoid $C_1$ is connected (also called transitive)) and $H_n(C, x) = 0$ for all $n \geqslant 2, x \in C_0$.

We assume the basic facts on crossed complexes as surveyed in for example [**Bro99, Bro04**]. In particular, we will use the monoidal closed structure on the category $\mathsf{Crs}$ which gives an exponential law of the form

$$\mathsf{Crs}(A \otimes B, C) \cong \mathsf{Crs}(A, \mathsf{CRS}(B, C)),$$

for crossed complexes $A, B, C$, and also in the pointed form

$$\mathsf{Crs}_*(A \otimes_* B, C) \cong \mathsf{Crs}_*(A, \mathsf{CRS}_*(B, C))$$

for pointed crossed complexes $A, B, C$; both of these structures are developed in [**BH87**]. Here the elements of $\mathsf{CRS}(B, C)_0$ are morphisms $B \to C$; the elements of the groupoid $\mathsf{CRS}(B, C)_1$ are (left) homotopies between these morphisms; and the further elements of $\mathsf{CRS}(B, C)$ are forms of higher homotopies.

A full exposition of the theory of crossed complexes will be given in [**BHS08**].

## 2. Fibrations of crossed complexes

A model category structure on the category $\mathsf{Crs}$ has been studied by Brown and Golasiński in [**BG89**], exploiting the notion of fibration of crossed complexes as defined in [**How79**]. We recall some of this material, but with a slightly different emphasis.

Recall that a morphism $p : E \to B$ of groupoids is a *fibration (covering morphism)* if it is star surjective (bijective) [**Bro06**]. The extension to crossed complexes is quite simple (covering morphisms of crossed complexes are applied in [**BRS99**]).

**Definition 2.1.** A morphism $p : E \to B$ of crossed complexes is a *fibration* if

(i) the morphism $p_1 : E_1 \to B_1$ is a *fibration* of groupoids;

(ii) for each $n \geqslant 2$ and $x \in E_0$, the morphism of groups $p_n : E_n(x) \to B_n(px)$ is surjective.

The morphism $p$ is a *trivial fibration* if it is a fibration, and also a weak equivalence, by which is meant that $p$ induces a bijection on $\pi_0$ and isomorphisms $\pi_1(E, x) \to \pi_1(B, px)$, $H_n(E, x) \to H_n(B, px)$ for all $x \in E_0$ and $n \geqslant 2$. □

We now follow model category ideas as in [**BG89**].

**Definition 2.2.** Consider the following diagram.

If given $i$ the dotted completion exists for all morphisms $p$ in a class $\mathcal{F}$, then we say that $i$ has the *left lifting property (LLP)* with respect to $\mathcal{F}$. We say a morphism $i : A \to C$ is a *cofibration* if it has the LLP with respect to all trivial fibrations. We say a crossed complex $C$ is *cofibrant* if the inclusion $\emptyset \to C$ is a cofibration. □

To give our most important example important example of a cofibration, we need some definitions.

**Definition 2.3.** We write $\mathbb{C}(n)$ for the free crossed complex on one generator $c^n$ of dimension $n$, and $\mathbb{S}^{n-1}$ for the subcomplex of $\mathbb{C}(n)$ generated by the elements of dimension $(n-1)$. Thus $\mathbb{C}(1)$ is essentially a groupoid, often written $\mathcal{I}$, which, with the inclusions $0, 1 \to \mathbb{C}(1)$, is a unit interval object in the category of crossed complexes, [**KP97**]; $\mathbb{C}(n)$ is a model of the $n$-disc, and $\mathbb{S}^{n-1}$ is a model of the $(n-1)$-sphere.

**Definition 2.4.** Let $A$ be a crossed complex. A morphism $i : A \to F$ of crossed complexes is said to be *relatively free* if $i$ is the canonical morphism when $F$ is obtained by attaching in order of increasing dimension copies of $\mathbb{C}(n)$ by means of morphisms $\mathbb{S}(n-1)$, analogously to the corresponding notion for relative $CW$-complexes. The images of the elements $c^n$ are called basis elements of $F$. In the case $A$ is empty, then $F$ is called a *free crossed complex*.

The following is [**BG89**, Corolllary 2.4]. We give the proof for the convenience of the reader.

**Proposition 2.5.** *Let $i : A \to F$ be a relatively free morphism of crossed complexes. Then $i$ is a cofibration.*

*Proof.* We consider the following diagram

$$
\begin{array}{ccc}
A & \xrightarrow{\ \alpha\ } & E \\
{\scriptstyle i}\downarrow & {\scriptstyle g}\nearrow & \downarrow{\scriptstyle p} \\
F & \xrightarrow[\ f\ ]{} & B
\end{array}
$$

in which $p$ is supposed a trivial fibration, and the morphisms $f, \alpha$ satisfy $fi = p\alpha$. We construct the regular completion $g$ on a relatively free basis $X$ of $F$ by induction.

For $n = 0$, we just lift a point in $B$ to a point in the corresponding component in $E$. This defines $g^0$ on $A \cup X^0$.

For the case $n = 1$, consider a basis element $x \in X^1(a, b)$, so that $f(x) \in B(fa, fb)$, and $g^0a, g^0b$ belong to the same component of $E$, by the condition on $\pi_0$. So there is an element $e \in E(g^0a, g^0b)$. Hence $fx - pe$ is a loop in $B_1(g^0a)$. By the condition for $p$ on $\pi_1$, there is a loop $e_1 \in E_1(g^0a)$ such that $pe_1$ is equivalent to $fx - pe$, i.e. $pe_1 = fx - pe + \delta_2 b_2$ for some $b_2 \in B_2(fa)$. By the fibration condition, $b_2 = pe_2$ for some $e_2 \in E_2(g^0a)$. Then $p(e_1 + \delta_2 e_2 + e) = fx$. So we can choose $g^1x = e_1 + \delta_2 e_2 + e$ to obtain an extension on $x$.

Suppose $n \geqslant 2$ and $g$ is defined on $X^{n-1}$. Consider an element $x$ of the free basis in dimension $n$. Then $g\delta x$ is defined and $pg\delta x = f\delta x$.

By the fibration condition, we can choose $e_n \in E_n$ such that $pe_n = fx$. Let $w = g\delta x - \delta e_n \in E_{n-1}$. Then $pw = 0, \delta w = 0$. By the triviality condition, $w$ is a boundary, i.e. $w = \delta z$ for some $z \in E_n$. Then $\delta(z + e_n) = g\delta x$. So we can extend $g$ by defining it on $x$ to be $z + e_n$. $\qquad\square$

The following is [**BG89**, Proposition 2.1.] with a different proof.

**Proposition 2.6.** *The following are equivalent for a morphism* $p : E \to B$ *in* Crs*:*

  (i)  *$p$ is a trivial fibration:*

  (ii)  *$p_0$ is surjective; if $e, e' \in E_0$ and $b_1 \in B_1(p_0 e, p_0 e')$, then there is $e_1 \in E_1(e, e')$ such that $p_1 e_1 = b_1$; if $n \geqslant 1$ and $e \in E_n$ satisfies $\delta^0 e = \delta^1 e$ for $n = 1$, $\delta e = 0$ for $n \geqslant 2$, and $b \in B_{n+l}$ satisfies $\delta b = p_n e$, then there is*

$$e' \in E_{n+1} \quad \text{such that} \quad p_{n+1} e' = b \quad \text{and} \quad \delta e' = e;$$

  (iii)  *$p$ has the RLP with respect to $\mathbb{S}(n-1) \to \mathbb{C}(n)$ for all $n \geqslant 0$;*

  (iv)  *if $F$ is a free crossed complex then $p$ has the RLP with respect to $\mathbb{S}(n-1) \otimes F \to \mathbb{C}(n) \otimes F$ for all $n \geqslant 0$;*

  (v)  *if $F$ is a free crossed complex then the induced morphism $p_* : \mathsf{CRS}(F, E) \to \mathsf{CRS}(F, B)$ is a trivial fibration.*

*Proof.* The proof that (i), (ii), (iii) are equivalent follows easily from the definition and Proposition 2.5.

Further (iii) implies (iv) since we know that under the condition of (iv),

$$\mathbb{S}(n-1) \otimes F \to \mathbb{C}(n) \otimes F$$

is relatively free [**BH91**, Proposition 5.1]. Finally (iv) trivially implies (iii).   □

We also need a pointed version of (iv) of the previous proposition.

**Proposition 2.7.** *If $F, E, B$ are pointed reduced crossed complexes with $F$ free, and $p : E \to B$ is a trivial fibration, then so also is*

$$p_* : \mathsf{CRS}_*(F, E) \to \mathsf{CRS}_*(F, B).$$

*Proof.* This relies on the pointed exponential law, [**BH87**], and the clear fact that $\mathbb{S}^{n-1} \otimes_* F \to \mathbb{C}(n) \otimes_* F$ is relatively free, which follows from methods analogous to those of [**BH91**].   □

**Definition 2.8.** If $G$ is a groupoid, we write $\mathbb{K}(G, 1)$ for the crossed complex which is $G$ in dimension 1 and trivial elsewhere. Thus $\mathbb{K}(G, 1)$ is certainly aspherical.   □

**Proposition 2.9.** *Let $F, C$ be reduced crossed complexes with $F$ free and $C$ aspherical. Let $G = \pi_1(C)$. Then there are bijections*

$$\pi_0 \mathsf{CRS}_*(F, C) = [F, C]_* \cong [F, \mathbb{K}(G, 1)]_* \cong \mathrm{Hom}(\pi_1 F, \pi_1 C),$$

*and for all $f : F \to C$ and $n \geqslant 1$ we have*

$$\pi_n(\mathsf{CRS}_*(F, C), f) = 0.$$

*Proof.* Let $G = \pi_1 C$. Since $C$ is aspherical, the natural morphism $p : C \to \mathbb{K}(G, 1)$ is not only a fibration but also a weak equivalence of crossed complexes. It follows that $p_* : \mathsf{CRS}_*(F, C) \to \mathsf{CRS}_*(F, \mathbb{K}(G, 1))$ is a trivial fibration and so a weak equivalence. In particular, $p_*$ induces a bijection of $\pi_0$. This gives the first result, since $\pi_0 \mathsf{CRS}_*(F, \mathbb{K}(G, 1))$ is clearly bijective with $\mathrm{Hom}(\pi_1 F, G)$.

The second result follows, since all homotopies and higher homotopies $F \to \mathbb{K}(G, 1)$ are trivial.   □

**Remark 2.10.** The classical homological algebra type of inductive proof is somehow hidden in this proof. □

Howie states in [**How79**] that a fibration $p : E \to B$ of crossed complexes yields a family of exact sequences involving the $H_n$, $\pi_1$ and $\pi_0$, akin to the well known family of homotopy exact sequences of a fibration of spaces, or of groupoids; for the latter see [**Bro06**, 7.2.9]. The standard properties of these sequences are:

- dependency on base points;

- non abelian features in dimension 1;

- sets with base point in dimension 0 but with some useful information obtainable from operations.

Let $x \in E_0$ and let $\mathcal{F}_x = p^{-1}(px)$ be the sub crossed complex of $E$ consisting of all elements of $E_0$ which map by $p$ to $x$ and of all elements of some $E_n, n \geqslant 1$, which map by $p$ to the identity at $px$. Here is Howie's exact sequence.

**Theorem 2.11.** *There is an exact sequence*

$$\cdots \to H_n(\mathcal{F}_x, x) \xrightarrow{i_n} H_n(E, x) \xrightarrow{p_n} H_n(B, px) \xrightarrow{\partial_n} \cdots$$

$$\cdots \to \pi_1(\mathcal{F}_x, x) \xrightarrow{i_1} \pi_1(E, x) \xrightarrow{p_1} \pi_1(B, px) \xrightarrow{\partial_1} \pi_0(\mathcal{F}_x) \xrightarrow{i_*} \pi_0(E) \xrightarrow{p_*} \pi_0(B).$$

*Here the terms of the sequence are all groups, except the last three which are sets with base points the classes $x_{\mathcal{F}}, x_E, x_B$ of $x$, $x$, $px$ respectively.*

(i) *There is an operation of the group $\pi_1(E, x)$ on the group $\pi_1(\mathcal{F}_x, x)$ making the morphism*

$$i_1 : \pi_1(\mathcal{F}_x, x) \to \pi_1(E, x)$$

*into a crossed module.*

(ii) *There is an operation of the group $\pi_1(B, px)$ on the set $\pi_0(\mathcal{F}_x)$ such that the boundary $\pi_1(B, px) \xrightarrow{\partial_1} \pi_0(\mathcal{F}_x)$ is given by $\partial_1(\alpha) = \alpha \cdot x_F$.*

*Further we have additional exactness at the bottom end as follows:*

(a) *$\partial_1 \alpha = \partial_1 \beta$ if and only if there is a $\gamma \in E(x)$ such that $p_1 \gamma = -\beta + \alpha$;*

(b) *if $\bar{u}$ denotes the component in $\mathcal{F}_x$ of an object $u$ of $\mathcal{F}_x$, then $i_* \bar{u} = i_* \bar{v}$ if and only if there is an $\alpha \in B(y)$ such that $\alpha_\# \bar{u} = \bar{v}$;*

(c) *if $\hat{y}$ denotes the component of $y$ in $B$ then*

$$i_*[\pi_0 \mathcal{F}_x] = p_*^{-1}[\hat{y}].$$

*Proof.* The proof of this theorem is a development of the part of the theorem which deals with fibrations of groupoids and which is given for example by Brown in [**Bro06**]. We leave the details as an exercise. □

**Remark 2.12.** It may be useful to point out that the exact sequence of a fibration of groupoids is generalised to a Mayer-Vietoris type sequence for a pullback of a fibration in [**BHK83**]. See also [**Bro06**, 10.7.6].

## 3.   Homotopy classification of maps to an $n$-aspherical space

We have already analysed a simple case of $[F, C]_*$ in Proposition 2.9. We now proceed to some slightly more complicated examples, using the exact sequences of a fibration. We also concentrate on the pointed homotopy classification, since this is nearest to classical results on group extensions, but the unpointed case is handled similarly.

**Definition 3.1.** A space $X$ is called *$n$-aspherical* if it is path-connected and $\pi_i(X) = 0$ for $1 < i < n$. Thus any path-connected space is 2-aspherical. Analogous terms are applied also to crossed complexes.

**Definition 3.2.** For a group or groupoid $Q$, $Q$-module $A$, and $n \geqslant 2$, let $\mathbb{K}(Q, 1; A, n)$ denote the crossed complex which is $Q$ in dimension 1, $A$ in dimension $n$, with the given action of $Q$, and all boundaries are trivial.     □

**Remark 3.3.** Let $F$ be a crossed complex. We will usually write $\phi : F_1 \to \pi_1 F$ for the quotient morphism. Then for any morphism $\theta : \pi_1 F \to Q$ the composite $\theta\phi$ is completely determined by $\theta$. However sometimes we start with a group (or groupoid) $\Phi$ and choose what we have called a *free crossed resolution $F$ of* $\Phi$. This is an aspherical free crossed complex together with a choice of isomorphism $\pi_1 F \to \Phi$, or, equivalently, with a quotient morphism $\phi : F_1 \to \Phi$ with kernel the image of $\delta_2 : F_2 \to F_1$. In such case $\phi$ is not determined by $F$.

**Proposition 3.4.** *Let $p : E \to B$ be a morphism of crossed complexes. Then $p$ is a fibration if and only if for any free crossed complex $F$, the induced morphism $p_* : \mathsf{CRS}(F, E) \to \mathsf{CRS}(F, B)$ is a fibration. If further $p$ and $F$ are pointed, then the induced morphism of pointed internal homs $p_* : \mathsf{CRS}_*(F, E) \to \mathsf{CRS}_*(F, B)$ is a fibration of crossed complexes.*

*Proof.* The forward implication follows easily from the exponential law, as in the proof of Proposition 2.6.

To prove it in the other direction, one takes again $F$ to be $\mathbb{C}(n)$, the free crossed complex on one generator of dimension $n$.     □

**Definition 3.5.** *Let $F, C$ be crossed complexes, let $A$ be a subcomplex of $F$, with inclusion $i : A \to F$ and let $f : A \to C$ be a morphism. We write $[F, C; f]$ for the set of homotopy classes rel $A$ of morphisms $F \to C$ which extend $f$. We write similarly $[F, C; f]_*$ for the pointed homotopy classes in the case $A, F, C, i$ are pointed.*

**Theorem 3.6.** *Let $F$ be a reduced free crossed complex and let $\Phi = \pi_1 F$. Then $[F, \mathbb{K}(Q, 1; A, n)]_*$ is the disjoint union of sets $[F, \mathbb{K}(Q, 1; A, n) : \theta\phi]_*$ one for each morphism $\theta : \Phi \to Q$, namely those homotopy classes inducing $\theta$. Further, the morphisms $F \to \mathbb{K}(Q, 1; A, n)$ inducing $\theta : \Phi \to Q$ may be given the structure of abelian group which is inherited by homotopy classes.*

*Proof.* The morphism $q : \mathbb{K}(Q, 1; A, n) \to \mathbb{K}(Q, 1)$ which is the identity in dimension 1 and 0 elsewhere is a fibration inducing a fibration

$$q_* : \mathsf{CRS}_*(F, \mathbb{K}(Q, 1; A, n)) \to \mathsf{CRS}_*(F, \mathbb{K}(Q, 1)).$$

The induced map on $\pi_0$ is surjective since every morphism $f : F \to \mathbb{K}(Q, 1)$ may be lifted by 0 to a morphism $F \to \mathbb{K}(Q, 1; A, n)$. By Proposition 2.9, $\pi_0 \mathsf{CRS}_*(F, \mathbb{K}(Q, 1)) \cong \mathrm{Hom}(\Phi, Q)$. So we can write, using the exact sequence of Theorem 2.11,

$$\pi_0 \mathsf{CRS}_*(F, \mathbb{K}(Q, 1; A, n)) \cong \bigsqcup_{\theta : \Phi \to Q} [F, \mathbb{K}(Q, 1; A, n); \theta\phi]_*.$$

The abelian group structure on each set $[F, \mathbb{K}(Q, 1; A, n); \theta\phi]_*$ by addition of values in dimension $n$ is clear from the diagram

**Definition 3.7.** We write $H^n_{\theta\phi}(F, A)$ for $[F, \mathbb{K}(Q, 1; A, n); \theta\phi]$, and call this abelian group the *nth cohomology over $\theta\phi$* of $F$ with coefficients in $A$. Thus $[F, \mathbb{K}(Q, 1; A, n)]_*$ is the disjoint union of the abelian groups $H^n_{\theta\phi}(F, A)$ for all morphisms $\theta : \Phi \to Q$. When convenient and clear, we abbreviate $\theta\phi$ to $\theta$.

A generalisation of the previous example is as follows.

**Example 3.8.** *Let $C$ be a reduced crossed complex such that $C_1 = Q$, and $\delta_2 = 0 : C_2 \to C_1$. Let $F$ be a free crossed complex. Then $\mathsf{Crs}_*(F, C)$ and $[F, C]_{\theta\phi}$ may be given the structure of abelian group by addition of values.* □

We now obtain for homotopy classification of maps a result which is analogous to and in fact directly generalises the classical theory of abstract kernels and obstructions, [**ML63**, Ch.IV, Thm.8.7]. First we give a definition.

**Definition 3.9.** Let $F$ be a free reduced crossed complex, and let $\phi : F_1 \to G = \pi_1(F)$ be the canonical morphism. If $\theta : G \to Q$ is a morphism of groups, and $A$ is a $Q$-module, and $n \geqslant 2$, we define the *nth cohomology of $F$ with coefficients in $A$ with respect to $\theta\phi$* to be the abelian group

$$H^n_{\theta\phi}(F, A) = [F, \mathbb{K}(Q, 1; A, n) : \theta\phi]_*.$$

**Theorem 3.10.** *Let $n \geqslant 2$ and let $F, C$ be reduced crossed complexes such that $F$ is free, $C$ is n-aspherical, and $C_i = 0$ for $i > n$. Let $\Phi = \pi_1(F), Q = \pi_1 C, A = \mathrm{Ker}\, \delta_n : C_n \to C_{n-1}$. Let $\theta : \Phi \to Q$ be a morphism of groups. Then there is defined an element $k_\theta \in H^{n+1}_{\theta\phi}(F, A)$, called the* obstruction class of $\theta$*, such that the vanishing of $k_\theta$ is necessary and sufficient for $\theta$ to be realised by a morphism $F \to C$.*

*If $k_\theta = 0$, then the set $[F, C; \theta\phi]$ of homotopy classes of morphisms $F \to C$ realising $\theta\phi$ is bijective with $H^n_{\theta\phi}(F, A)$.*

*Proof.* Consider the morphisms of crossed complexes $C \xrightarrow{j} \xi C \xrightarrow{p} \zeta C$ as shown in the following diagram:

Then $\xi C$ is aspherical, $\zeta C = \mathbb{K}(Q, 1; A, n+1)$, and $p : \xi C \to \zeta C$ is a fibration of crossed complexes.

Since $F$ is a free reduced crossed complex, we have an induced fibration of crossed complexes

$$p_* : \mathsf{CRS}_*(F, \xi C) \to \mathsf{CRS}_*(F, \zeta C). \tag{2}$$

On applying $\pi_0$ to this we get, considering previous identifications, a map of sets

$$p_* : \mathrm{Hom}(\Phi, Q) \longrightarrow \bigsqcup_{\theta \in \mathrm{Hom}(\Phi, Q)} H^{n+1}_{\theta\phi}(F, A). \tag{3}$$

**Lemma 3.11.** *A group morphism $\theta : \Phi \to Q$ maps to 0 in $H^{n+1}_{\theta\phi}(F, A)$ if and only if $\theta$ is induced by a morphism $F \to C$.*

*Proof.* Suppose $\theta$ is induced by a morphism $f : F \to C$. Then $f$ factors through $pj$ and is therefore 0 in $H^{n+1}_{\theta\phi}(F, A)$.

Suppose conversely that $\theta$ determines 0 in $H^{n+1}_{\theta\phi}(F, A)$. We know that $\theta$ is induced by a morphism $f' : F \to \xi C$. Then $pf$ is homotopic to 0 and so by the fibration condition $f'$ is homotopic to $f''$ such that $pf'' = 0$. Hence $f''$determines $f : F \to C$ such that $jf = f''$. Then $f$ also induces $\theta$. This proves the lemma. □

Let $\mathcal{F}(f)$ denote the fibre of $p_*$ over $pf$. Then we have an exact sequence

$$\to \pi_1(\mathsf{CRS}_*(F, \xi C), f) \to \pi_1(\mathsf{CRS}_*(F, \zeta C), pf)$$
$$\to \pi_0\mathcal{F}(f) \to \pi_0\mathsf{CRS}_*(F, \xi C) \to \pi_0\mathsf{CRS}_*(F, \zeta C).$$

By Proposition 2.9, $\pi_1(\mathsf{CRS}_*(F, \xi C), f) = 0$, and so the above sequence translates to

$$0 \to H^n_{\theta\phi}(F, A) \to [F, C; \theta\phi] \to \mathrm{Hom}(\Phi, Q).$$

Further we have a free action of the abelian group $H^n_{\theta\phi}(F, A)$ on the set $[F, C]_{\theta\phi}$. This completes the proof of the theorem. □

This result generalises the classical theory of extensions of groups and $Q$-kernels. To apply the theory to that case, the crossed complex $F$ is taken to be a free crossed resolution of the group $G$. If $F$ is the standard free crossed resolution of $G$, then the relation with factor systems is shown in [**BP96**]. The advantage of this approach is that it is clear that the standard free crossed resolution may be replaced by any free crossed resolution of $G$, and in many cases it is possible to construct small such resolutions; in these cases, for example if $G$ is $FP_3$, we may obtain a finite description of the classes of extensions.

## 4. Applications to spaces

For the applications to spaces we need to know when a space is of the homotopy type or homotopy $n$-type of $BC$ for some crossed complex $C$. It is shown in [**BH81**] that for any crossed complex $C$ there is a filtered space $Y_*$ such that $\Pi Y_* \cong C$. Whitehead in [**Whi49**] gives an example of a 5-dimensional free crossed complex which is not isomorphic to $\Pi Y_*$ for any $CW$-filtration $Y_*$.

A key result is the following:

**Theorem 4.1.** *Let $Y_*$ be the skeletal filtration of a $CW$-complex $Y$. Then there is 1-equivalence $q : Y \to B\Pi Y_*$ such that if $y \in Y_0$ then the exact sequence of the homotopy fibre over $y$ is equivalent to Whitehead's exact sequence [**Whi50**]:*

$$\cdots \to \Gamma_n(Y, y) \to \pi_n(Y, y) \xrightarrow{\omega} H_n(\widetilde{Y}_y) \to$$
$$\cdots \to \Gamma(\pi_2(Y, y)) \to \pi_3(Y, y) \xrightarrow{\omega} H_3(\widetilde{Y}_y) \to 0$$

*where $\omega$ is the Hurewicz morphism.*

*Proof.* The original form of this comes from [**BH81**, section 8], which in essence used a cubical version of the classifying space. The simplicial version is dealt with in [**Ash88**]. $\blacksquare$

**Corollary 4.2.** *If further $Y$ is $n$-aspherical, then the homotopy fibre of $q : Y \to B\Pi Y_*$ is $n$-connected.*

**Remark 4.3.** Whitehead in [**Whi49**] obtains his classification results on maps $X \to Y$ where $X$ is $n$-dimensional and when $Y$ is $J_n$-complex, which is essentially the condition that $q$ above is an $n$-equivalence, and which generalises the condition of $n$-asphericity. $\square$

Ronald Brown   `r.brown@bangor.ac.uk`

School of Computer Science, Bangor University, Dean St., Bangor, LL57 1UT, UK.